\documentclass[11pt]{article}

\usepackage[linesnumbered,ruled]{algorithm2e}

\usepackage{jhtitle}
\usepackage{geometry}
\usepackage{amsbsy,amsmath,amsthm,amssymb,graphicx}
\usepackage{natbib}
\usepackage{hyperref}
\usepackage{enumitem}
\usepackage{pstricks}
\usepackage{subcaption}

\newcommand{\pcite}[1]{\citeauthor{#1}'s \citeyearpar{#1}}

\newcommand{\df}{\mathrm{d}}

\newtheorem{theorem}{Theorem}

\newtheorem{remark}[theorem]{Remark}
\newtheorem{proposition}[theorem]{Proposition}

\newcommand{\bbe}{{\boldsymbol{\beta}}}

\title{\bf Approximating the Spectral Gap of the \\ {P}\'olya-Gamma
  Gibbs Sampler}

\author{
  Bryant Davis \\ Department of Statistics \\
  University of Florida \\
  \texttt{davibf11@ufl.edu}
  \and
  James P. Hobert \\ Department of Statistics
  \\ University of Florida \\
  \texttt{jhobert@stat.ufl.edu}
}

\date{April 2021}

\begin{document}
\maketitle

\begin{abstract}
  The self-adjoint, positive Markov operator defined by the
  {P}\'olya-Gamma Gibbs sampler (under a proper normal prior) is shown
  to be trace-class, which implies that all non-zero elements of its
  spectrum are eigenvalues.  Consequently, the spectral gap is
  $1-\lambda_*$, where $\lambda_* \in [0,1)$ is the second largest
    eigenvalue.  A method of constructing an asymptotically valid
    confidence interval for an upper bound on $\lambda_*$ is developed
    by adapting the classical Monte Carlo technique of
    \citet{qin:hobe:khar:2019} to the {P}\'olya-Gamma Gibbs sampler.
    The results are illustrated using the German credit data.  It is
    also shown that, in general, uniform ergodicity does not imply the
    trace-class property, nor does the trace-class property imply
    uniform ergodicity.
\end{abstract}

\section{Introduction}
\label{sec:intro}

Let $Y_1, \ldots, Y_n$ be independent Bernoulli random variables such
that $P(Y_i = 1) = F(x_i^T \beta)$ where $x_i$ is a $p \times 1$
vector of known covariates associated with $Y_i$, $\beta$ is an
unknown $p \times 1$ vector of regression parameters, and $F(x) =
e^x/(1+e^x)$, which is the standard logistic distribution function.
Under this model, the joint mass function of $Y_1, \ldots, Y_n$ is
\begin{align*}
    \prod_{i=1}^{n} P(Y_i = y_i \mid \beta) = \prod_{i=1}^{n} \big[
      F(x_i^{T} \beta) \big]^{y_i} \big[ 1 - F(x_i^{T} \beta) \big]^{1
      - y_i} I_{\{0, 1 \}}(y_i) \,.
\end{align*}
Let $y$ denote the $n \times 1$ vector of observed $y_i$ values.  We
consider a Bayesian analysis of this data under a proper prior for the
unknown vector $\beta$.  In particular, we take the prior on $\beta$
to be $\mbox{N}_p(b, B)$, where $b \in \mathbb{R}^p$ and $B$ is a $p
\times p$ positive definite matrix.  The resulting posterior density,
$\pi(\beta \,|\, y)$, is highly intractable, but there is a simple
MCMC algorithm that can be used to explore it.  This algorithm was
developed by \citet{pols:scot:wind:2013}, and is often called the
{P}\'olya-Gamma Gibbs sampler.  In order to describe the algorithm, we
must introduce some notation.

Let $X$ be the $n \times p$ matrix whose $i$th row is $x_i^T$, and
define $\mathbb{R}_+ = (0,\infty)$.  For $w \in \mathbb{R}_+^n$, let
$\Omega(w) = \mbox{diag} \big \{ w_i \big\}_{i=1}^n$.  We then define
\[
\Sigma(w) = \left( X^T \Omega(w) X + B^{-1} \right)^{-1} \;,
\]
as well as
\[
\mu(w) = \Sigma(w) \bigg[ X^T \bigg( y - \frac{1}{2} 1_n \bigg) +
  B^{-1} b \bigg] \;,
\]
where $1_n$ is an $n \times 1$ vector of 1s.  When we write $Z \sim
\mbox{PG}(1, d)$, we mean that the random variable $Z$ follows a
P\'olya-Gamma distribution with probability density function (pdf)
given by
\[
f(z ; d) = \mbox{cosh}(d/2) e^{-d^2 z/2} g(z) \;,
\]
where $d \ge 0$, and
\[
g(z) = \sum_{k=0}^{\infty} (-1)^{k} \; \frac{2k + 1}{\sqrt{2 \pi z^3}}
\; \mbox{exp} \left\{ -\frac{(2k + 1)^2}{8z} \right\} I_{(0,
  \infty)}(z) \;.
\]
The function $g(z)$ is a pdf \citep{bian:pitm:yor:2001}, and, in
particular, it is the pdf of a $\mbox{PG}(1, 0)$ random variable (for
more on the P\'olya-Gamma distribution, including highly efficient
methods for simulating from it, see \citet{pols:scot:wind:2013} and
\citet{wind:pols:scot:2014}).  The {P}\'olya-Gamma (PG) Gibbs sampler
simulates a Markov chain $\Gamma = \{\beta^{(m)}\}_{m=0}^\infty$ using
the following two-step procedure to move from the current state,
$\beta^{(m)}$, to the new state, $\beta^{(m+1)}$.

\begin{algorithm}
    \vspace*{2.5mm}
    \begin{itemize}
        \item[1.] Draw $W_1, W_2, \ldots, W_n$ independently with 
        $$ W_i \sim \mbox{PG}(1, \vert x_i^T \beta^{(m)} \vert) \,, $$
          and call the observed vector $w = (w_1, \ldots, w_n)^T$.
        \item[2.] Draw $\beta^{(m+1)} \sim \mbox{N}_{p} \big( \mu(w),
          \Sigma(w) \big)$ \,.
    \end{itemize}
    \caption{Iteration $m + 1$ of the PG Gibbs sampler.}
\end{algorithm}

The Markov chain $\Gamma$ is irreducible, aperiodic and positive
Harris recurrent, and the posterior density, $\pi(\beta \,|\, y)$, is
its unique invariant density \citep[see, e.g.,][]{choi:hobe:2013}.  In
this paper, we study the self-adjoint, positive Markov operator
associated with $\Gamma$, call it $K$.  Our main theoretical result is
that this operator is \textit{trace-class} (for background on the
spectral theory of linear operators, see, e.g., \citet{helm:2014} and
\citet{conw:1990}, and for Markov operators in particular, see
\citet{mira:geye:1999} and \citet{qin:hobe:khar:2019}).  The
trace-class property implies that all non-zero elements of the
spectrum of $K$ are eigenvalues, so the spectral gap of $K$, which
controls the rate at which $\Gamma$ converges to $\pi(\beta \,|\, y)$,
is equal to $1-\lambda_*$, where $\lambda_* \in [0,1)$ is the second
  largest eigenvalue.  Furthermore, the trace-class property also
  implies that the sum of the eigenvalues is finite, and this fact
  allows us to use the results of \citet{qin:hobe:khar:2019} to
  develop a method for constructing an asymptotically valid confidence
  interval for an upper bound on $\lambda_*$.  Our point estimator of
  the upper bound on $\lambda_*$ is based on iid random vectors that
  can be simulated using non-stationary, short runs of the Markov
  chain $\Gamma$.  We illustrate the application of our method using a
  medium-sized data set concerning 1,000 applications for credit in
  Germany, the so-called German credit data.

Previous work on the theoretical properties of PG Gibbs samplers
includes \citet{choi:hobe:2013}, \citet{choi:roma:2017}, and
\citet{wang:roy:2018a}.  \citet{choi:hobe:2013} analyzed the PG Gibbs
sampler under the same normal prior on $\beta$ and showed that it is
uniformly ergodic.  It is well-known that the trace-class property
implies geometric ergodicity (see Section~\ref{sec:properties}), and
that uniform ergodicity implies geometric ergodicity; however, prior
to our work, the relationship between uniform ergodicity and the
trace-class property (for ergodic Markov chains with self-adjoint and
positive Markov operators) was unknown.  In
Section~\ref{sec:relationship}, we show that, in general, uniform
ergodicity does not imply the trace-class property, nor does the
trace-class property imply uniform ergodicity.  \citet{choi:roma:2017}
and \citet{wang:roy:2018a} analyzed a slightly different PG Gibbs
sampler.  In particular, instead of the proper normal prior on $\beta$
that we consider here, these authors used a flat improper prior.
\citet{wang:roy:2018a} showed that this version of the PG Gibbs
sampler is geometrically ergodic, and \citet{choi:roma:2017} showed
that, if each row of the $X$ matrix is a unit vector, then the
corresponding Markov operator is trace-class.  We note that when the
rows of $X$ are unit vectors (which corresponds to a one-way ANOVA
design) and the prior on $\beta$ is flat, MCMC is actually not
required because the $p$ (univariate) components of $\beta$ are
\textit{a posteriori} independent.  Indeed, the posterior density
factors into a constant times a product of $p$ terms, each having the
form $\big[ F(\beta_i) \big]^{a_i} \big[ 1 - F(\beta_i) \big]^{b_i}$,
where $a_i$ and $b_i$ are non-negative integers.  Hence, one could
presumably design a univariate accept-reject algorithm that would
yield exact draws from the posterior.  Finally, for extensions to
\textit{mixed} logistic regression models, see \citet{wang:roy:2018b}
and \citet{rao:roy:2021}.

The remainder of the paper is organized as follows.
Section~\ref{sec:properties} describes some theoretical properties of
the PG Gibbs sampler, and contains a statement of our main theoretical
result.  Our implementation of \pcite{qin:hobe:khar:2019} algorithm in
the context of the PG Gibbs sampler is described in
Section~\ref{sec:aprox}.  In Section~\ref{sec:example}, we use our
results to approximate the spectral gap of the PG Gibbs sampler for
the German credit data.  The relationship between uniform ergodicity
and the trace class property is the topic of
Section~\ref{sec:relationship}.  Two proofs are relegated to an
Appendix.

\section{Theoretical Properties of the PG Gibbs Sampler}
\label{sec:properties}

The Markov transition density (Mtd) of the PG Gibbs sampler can be
expressed as
\[
k(\beta,\beta') = \int_{\mathbb{R}_+^n} \pi_2(\beta' \,|\, w, y)
\pi_1(w \,|\, \beta, y) \, dw \;,
\]
where, using notation from the Introduction, $\pi_1(w \,|\, \beta, y)
= \prod_{i=1}^n f \big( w_i ; | x_i^{T} \beta | \big)$, and $\pi_2$ is
a $p$-variate normal density with mean $\mu(w)$ and variance
$\Sigma(w)$.  Of course, $\pi_1$ and $\pi_2$ are the conditional pdfs
corresponding to an \textit{augmented} posterior density,
$\pi_a(\beta, w \,|\, y)$, which satisfies $\int_{\mathbb{R}_+^n}
\pi_a(\beta, w \,|\, y) \, dw = \pi(\beta \,|\, y)$.  (See
\citet{choi:hobe:2013} for the exact form of $\pi_a(\beta, w \,|\,
y)$.)  We can also define the \textit{conjugate} Markov chain,
$\tilde{\Gamma} = \{w^{(m)}\}_{m=0}^\infty$, which lives on
$\mathbb{R}_+^n$, and has Mtd given by
\[
\tilde{k}(w,w') = \int_{\mathbb{R}^p} \pi_1(w' \,|\, \beta, y)
\pi_2(\beta \,|\, w, y) \, d\beta \;.
\]
The chain $\tilde{\Gamma}$ is also irreducible, aperiodic and positive
Harris recurrent, and its unique invariant density is
$\int_{\mathbb{R}^p} \pi_a(\beta, w \,|\, y) \, d\beta$.  We will make
use of $\tilde{\Gamma}$ in the sequel.

We now turn to the convergence properties of the PG Gibbs sampler.
Let $L^2(\pi)$ denote the Hilbert space of complex valued functions
that are square integrable with respect to the target posterior
density, $\pi(\beta \,|\, y)$, i.e.,
\[
L^2(\pi) := \Big\{ f: \mathbb{R}^p \to \mathbb{C} \; \Big \arrowvert
\, \int_{\mathbb{R}^p} |f(\beta)|^2 \pi(\beta \,|\, y) \, d\beta <
\infty \Big\} \,.
\]
The inner product of $f, g \in L^2(\pi)$ is given by
\[
\langle f, g \rangle_{\pi} = \int_{\mathbb{R}^p} f(\beta)
\overline{g(\beta)} \pi(\beta \,|\, y) \, d\beta \,.
\]
The Mtd $k(\beta,\beta')$ defines a linear (Markov) operator~$K :
L^2(\pi) \rightarrow L^2(\pi)$ such that if $f \in L^2(\pi)$, then
\[
Kf(\beta) = \int_{\mathbb{R}^p} k(\beta,\beta') f(\beta') \, d\beta'
\,.
\]
The Markov chain $\Gamma$ is the $\beta$-marginal of the two-block
Gibbs chain (based on $\pi_a$) that alternates between $\beta$ and
$w$.  It follows that $K$ is positive (and self-adjoint)
\citep{liu:wong:kong:1994}.  Our main theoretical result, which is
proven in Appendix~\ref{app:A}, shows that $K$ enjoys additional
regularity.

\begin{proposition}
  \label{prop:tc}
  The Markov operator $K$ is trace-class.
\end{proposition}

The fact that $K$ is positive (and self-adjoint) implies that its
spectrum is contained in $[0,1]$.  The trace-class property implies
that $K$ is a compact operator, so it possesses a pure eigenvalue
spectrum with (at most) a countable number of eigenvalues.  Let
$\{\lambda_i\}_{i=0}^\kappa$ denote the strictly positive eigenvalues
of $K$ in decreasing order, where $0 \le \kappa \le \infty$.  We know
that $\lambda_0 = 1$ and, since there cannot be two eigenvalues equal
to 1 (because of irreducibility), we have $\lambda_1 \in [0,1)$.  (In
  what follows, we use the symbols $\lambda_1$ and $\lambda_*$
  interchangeably.)  Hence, the spectral gap is $1-\lambda_*>0$, so
  $\Gamma$ is geometrically ergodic and $\lambda_*$ can be viewed as
  its rate of convergence \citep{robe:rose:1997}.  Indeed, let $\Pi$
  denote the probability measure defined by $\pi(\beta \,|\, y)$, and
  let $\nu$ be any probability measure that is absolutely continuous
  with respect to $\Pi$ and satisfies $\int_{\mathbb{R}^p} (\df
  \nu/\df \Pi)^2 \, \df \Pi < \infty$.  Then there exists a constant
  $M_{\nu} < \infty$ such that
\begin{equation*}
  d_{\mbox{\scriptsize{TV}}}(\nu K^m, \Pi) \leq M_{\nu} \, \lambda_*^m
  \,,
\end{equation*}
where $\nu K^m$ denotes the probability measure of $\beta^{(m)}$ when
$\beta^{(0)} \sim \nu$, and $d_{\mbox{\scriptsize{TV}}}(\cdot,\cdot)$
is the total variation distance.  (As usual, the symbol $K$ is doing
double duty, representing both the Markov operator and the Markov
transition kernel.)

The trace-class property also implies that the eigenvalues of $K$ are
summable, and this means that we can employ the method of
\citet{qin:hobe:khar:2019} (hereafter QH\&K) to estimate $\lambda_*$.
The details are presented in the next section.

\section{Approximating $\lambda_*$}
\label{sec:aprox}

We begin with a high-level description of QH\&K's method.  As stated
above, $\sum_{i=0}^{\kappa} \lambda_i < \infty$ is a result of the
trace-class property of $K$, and it follows that, if $l$ is any
strictly positive integer, then $s_l := \sum_{i=0}^\kappa \lambda_i^l
< \infty$ as well.  QH\&K develop a classical Monte Carlo estimator of
$s_l$.  Moreover, they show that $u_l := (s_l - 1)^{1/l} \downarrow
\lambda_*$ as $l \to \infty$.  Thus, for fixed $l$, the estimator of
$s_l$ can be converted directly into an estimator of an upper bound on
$\lambda_*$.  QH\&K provide guidance on the choice of $l$, as well as
sufficient conditions under which their Monte Carlo estimator of $s_l$
has finite variance, which allows one to calculate a standard error
for $u_l$ via the delta method.

In order to use QH\&K's method, an auxiliary density, $h: \mathbb{R}^p
\rightarrow (0,\infty)$, must be specified.  This density plays a role
similar to that of the importance density in an importance sampling
algorithm.  The classical Monte Carlo (unbiased) estimator of $s_l$ is
given by
\begin{equation}
  \label{eq:mce}
  \hat{s}_l = \frac{1}{N} \sum_{i=1}^N \frac{\pi(\beta_i^* \mid
    w_i^*)}{h(\beta_i^*)} \;,
\end{equation}
where $N$ is the Monte Carlo sample size, and the random vectors
$\{(\beta_i^*,w_i^*)\}_{i=1}^N$ are iid and generated according to
Algorithm 2.

\begin{algorithm}
    \vspace*{2.5mm}
    \begin{itemize}
        \item[1.] Draw $\beta^* \sim h(\cdot)$. \\
        \item[2.] Draw $W_1, W_2, \dots, W_n$ independently with
        \[
        W_i \sim \mbox{PG}(1, \mid x_i^T \beta^* \mid) \,,
        \]
        and call the observed vector $w = (w_1,w_2,\dots,w_n)^T$. \\
        \item[3.] If $l = 1$, set $w^* = w$.  If $l \ge 2$, draw $w^*
          \sim \tilde{k}^{(l-1)}(\cdot \mid w)$ by running $l-1$ iterations
          of the conjugate chain $\tilde{\Gamma}$ initiated at $w$.
    \end{itemize}
    \caption{Drawing $(\bbe^*,w^*)$}
\end{algorithm}

The estimator \eqref{eq:mce} is strongly consistent no matter what $h$
is used.  However, additional conditions are required to guarantee
finite variance, which ensures the existence of asymptotically valid
confidence intervals for $s_l$ and $u_l$.  In particular, QH\&K show
that the following condition is sufficient for the estimator
\eqref{eq:mce} to have finite variance:
\[
  \int_{\mathbb{R}_+^n} \int_{\mathbb{R}^p} \frac{ \pi(w \mid \beta,
    y) \, \pi^3(\beta \mid w, y)}{h^2(\beta)} \; d\beta \; dw < \infty
  \;.
\]
The next result, which is proven in Appendix~\ref{app:B}, shows that
this condition is satisfied if $h$ is taken to be a multivariate
Student's $t$ density.

\begin{proposition}
  \label{prop:fv}
  Let $h_{\nu}(\beta ; d, C)$ denote a $p$-dimensional Student's $t$
  density with location parameter $d$, positive definite scale matrix
  $C$, and degrees of freedom $\nu$.  Then
\[
  \int_{\mathbb{R}_+^n} \int_{\mathbb{R}^p} \frac{ \pi(w \mid \beta,
    y) \, \pi^3(\beta \mid w, y)}{h_{\nu}^2(\beta ; d, C)} \; d\beta
  \; dw < \infty \;.
\]
\end{proposition}

\begin{remark}
  QH\&K actually present two different Monte Carlo estimators for
  $s_l$ - one is more effective for data sets in which $n$ is larger
  than $p$, while the other tends to work better when $p$ is larger
  than $n$.  In this paper, we consider only the former.
\end{remark}

In the next section, we use the results above to approximate the
spectral gap of the PG Gibbs sampler for a medium-sized real data set.

\section{An Application: The German Credit Data}
\label{sec:example}

In this section, we apply our method to the so-called German credit
data, which are available here:
\verb|http://archive.ics.uci.edu/ml/index.php|.  In this data set,
there are $n=1000$ binary observations, each one representing the
success or failure of a particular loan application.  (700 of the 1000
applications were deemed creditworthy.)  Associated with each of the
1,000 observations are twenty covariates including loan purpose,
demographic information, bank account balances, marital status, and
employment status.  Seven of the covariates are quantitative and
thirteen are categorical.  After converting the categorical covariates
to indicators, there are a total of $p=49$ regression coefficients.
This data set is frequently used to illustrate newly developed
statistical and machine learning techniques for binary data
\citep[see, e.g.,][]{pols:scot:wind:2013,jaco:olea:atch:2019}.

Our Bayesian model contains the hyper-parameters $b$ and $B$, which
must be specified.  Following \citet{jaco:olea:atch:2019}, we place a
relatively uninformative prior on $\beta$ with $b=0$ and $B = 10
I_{49}$.  In order to use Algorithm 2, the parameters of the Student's
$t$ density must be specified, and we used a preliminary run of the PG
Gibbs sampler to choose $d$ and $C$.  In particular, we ran the
sampler for 25,000 iterations (with $\beta^{(0)}$ set equal to the MLE
of $\beta$ based on the frequentist version of our logistic regression
model), discarded the first 5,000 draws as burn-in, and then used the
remaining 20,000 draws to get estimates (the usual ergodic averages)
of the posterior mean and covariance matrix, call these $\hat{\beta}$
and $\hat{\Sigma}$.  We then set $d = \hat{\beta}$, $C = \hat{\Sigma}$
and $\nu=5$.  Based on the guidance given in QH\&K, and some initial
experimentation, it became clear that $l=5$ was a reasonable choice.
We utilized a Monte Carlo sample size of $N = 10^7$.  The simulations
yielded $\hat{u}_5 = 0.787$ with a standard error of 0.074, resulting
in an asymptotically valid 95\% confidence interval of (0.639,0.935)
for $u_5$.  Hence, we can be fairly confident that the unknown
spectral gap $1-\lambda_*$ is at least $0.065$.

The next section concerns the relationship between uniform ergodicity
and the trace-class property.

\section{Uniform Ergodicity Versus the Trace-class Property}
\label{sec:relationship}

Consider an ergodic Markov chain whose Markov operator is self-adjoint
and positive.  As we explained in Section~\ref{sec:properties}, if
this Markov chain is trace-class, then it is also geometrically
ergodic.  A natural question that arose during our study of the
P\'{o}lya-Gamma Gibbs sampler is this: Is there a similar relationship
between uniform ergodicity and the trace-class property?  In this
section, we show that the answer is ``no.''  In particular, we show
that there exist chains that are uniformly ergodic, but not
trace-class, as well as chains that are trace-class, but not uniformly
ergodic.

Suppose we wish to sample from the univariate target density
\[
\pi(x) = \frac{2}{3}(x+1) I_{(0,1)}(x) \;.
\]
Consider an independence Metropolis algorithm based on a
$\mbox{Uniform}(0,1)$ candidate. The corresponding Markov operator is
clearly self-adjoint, and Lemma 3.1 of \citet{rudo:ullr:2013} implies
that it is also positive.  Now letting $q(y)$ denote the
$\mbox{Uniform}(0,1)$ candidate density, we can observe that for all
$y \in (0,1)$, we have
\[
\frac{q(y)}{\pi(y)} = \frac{3}{2(x+1)} \ge \frac{3}{4} \;.
\]
Hence, Theorem 2.1 of \citet{meng:twee:1996} implies that the chain is
uniformly ergodic.  However, the operator cannot be compact
\citep[see][p. 1755]{chan:geye:1994}, and thus cannot be trace-class.
Hence, we have a uniformly ergodic chain that is not trace-class.

For the other direction, we turn to a collection of birth-death Markov
chains that was analyzed in \citet{tan:jone:hobe:2013}.  Let
$\{a_i\}_{i=1}^\infty$ and $\{b_i\}_{i=1}^\infty$ be two sequences of
strictly positive real numbers satisfying
\[
\sum_{i=1}^\infty a_i + \sum_{i=1}^\infty b_i = 1 \;.
\]
Additionally, define $b_0 = 0$ and $\mathbb{N} = \{1,2,3,\dots\}$.  We
use the two sequences to define a bivariate random vector $(X,Y)$ with
support $\mathbb{N} \times \mathbb{N}$ and probability mass function
given by
\[ 
\pi(x, y) =
   \begin{cases} 
      a_x & x = y, y = 1, 2, 3, \ldots \\
      b_y & x = y+1, y = 1, 2, 3, \ldots \\
      0 & \mbox{otherwise} \;. 
   \end{cases}
\]
The marginal mass functions are $\pi_{X}(x) = a_x + b_{x-1}$ and
$\pi_{Y}(y) = a_y + b_y$, and the conditional mass functions are given
by
\[
 \pi_{X \,|\, Y}(x \,|\, y) = \frac{a_y}{a_y + b_y} I(x = y) +
 \frac{b_y}{a_y + b_y} I(x = y+1)
\]
for $y \in \mathbb{N}$, and 
\[
\pi_{Y \,|\, X}(y \,|\, x) = \frac{a_x}{a_x + b_{x-1}} I(y = x) +
\frac{b_{x-1}}{a_x + b_{x-1}} I(y=x-1)
\]
for $x \in \mathbb{N}$.  Let $\Gamma = \{X_n\}_{n=0}^\infty$ be the
Markov chain on $\mathbb{N}$ with Markov transition probabilities
given by
\[
P(X_{n+1} = x' | X_n = x) = k(x,x') = \sum_{y \in \mathbb{N}} \pi_{X
  \,|\, Y}(x' \,|\, y) \pi_{Y \,|\, X}(y \,|\, x) \;.
\]
Because this Markov chain is the marginal of a two component Gibbs
sampler, the corresponding operator is necessarily self-adjoint and
positive \citep{liu:wong:kong:1994}.  As we now explain, $\Gamma$
turns out to be a birth-death chain.  For $x \in \mathbb{N}$, define
\[
p_x = \frac{a_x b_x}{(a_x + b_{x-1})(a_x + b_x)} \;,
\]
and, for $x \in \{2,3,4,\dots\}$, define
\[
q_x = \frac{a_{x-1} b_{x-1}}{(a_x + b_{x-1})(a_{x-1} + b_{x-1})} \;.
\]
Finally, for $x \in \{2,3,4,\dots\}$, define $r_x = 1 - p_x - q_x$.
Using this notation, we can express $k$ as follows:
\[ 
k(x,x') = 
   \begin{cases} 
      1 - p_1 & x' = x = 1 \\
      p_x & x' = x + 1, \; x = 1, 2, 3, \ldots \\
      q_x & x' = x - 1, \; x = 2, 3, 4, \ldots \\
      r_x & x' = x, \hspace{0.83cm} x = 2, 3, 4, \ldots \\
      0 & \mbox{otherwise} \;. 
   \end{cases}
\]
Hence, at each step, the chain either goes up one unit, down one unit,
or stays at the current state.  Now, according to Theorem 2 in
\citet{qin:hobe:khar:2019}, $\Gamma$ is trace-class if and only if
\[
\sum_{x = 1}^{\infty} k(x, x) = 1 - p_1 + \sum_{x = 2}^{\infty} r_x <
\infty \;.
\]
The following representation of $r_x$ will be used below:
\begin{align*}
    r_x & = \frac{a_x^2}{(a_x + b_{x-1})(a_x + b_x)} +
    \frac{b_{x-1}^2}{(a_x + b_{x-1})(a_{x-1} + b_{x-1})} \\ & =
    \frac{1}{\left( 1 + \frac{b_{x-1}}{a_x} \right)\left( 1 +
      \frac{b_{x}}{a_x} \right)} + \frac{1}{\left( 1 +
      \frac{a_x}{b_{x-1}} \right)\left( 1 + \frac{a_{x-1}}{b_{x-1}}
      \right)} \;.
\end{align*}
We now consider specific versions of the sequences
$\{a_i\}_{i=1}^\infty$ and $\{b_i\}_{i=1}^\infty$, which are similar
to the ones employed in Example 2.5 of \citet{trun:2016}.  For $x \in
\mathbb{N}$, take $a'_x = (2x - 1)^{-(4x - 2)}$ and $b'_x =
(2x)^{-4x}$.  Let $c = \sum_{i=1}^\infty a'_i + \sum_{i=1}^\infty
b'_i$, and define $a^*_x = a'_x/c$ and $b^*_x = b'_x/c$.  Then, as
required, we have $\sum_{i=1}^\infty a^*_i + \sum_{i=1}^\infty b^*_i =
1$.  For $x \in \{2,3,4,\dots\}$, we have
\begin{align*}
    \frac{b^{*}_x}{a^{*}_x} = \frac{(2x - 1)^{4x - 2}}{(2x)^{4x}} =
    \left[ \frac{2x - 1}{2x} \right]^{4x - 2} \frac{1}{(2x)^2} \le
    \frac{1}{4x^2} \;,
\end{align*}
and
\begin{align*}
    \frac{a^{*}_{x+1}}{b^{*}_{x}} = \frac{(2x)^{4x}}{(2(x+1) -
      1)^{4(x+1) - 2}} = \left[ \frac{2x}{2x + 1} \right]^{4x}
    \frac{1}{(2x + 1)^2} \le \frac{1}{4x^2} \;.
\end{align*}
Hence, for $x \in \{2,3,4,\dots\}$, we see that
\begin{align*}
    r_x & = \frac{1}{\left( 1 + \frac{b^{*}_{x-1}}{a^{*}_x}
      \right)\left( 1 + \frac{b^{*}_{x}}{a^{*}_x} \right)} +
    \frac{1}{\left( 1 + \frac{a^{*}_x}{b^{*}_{x-1}} \right)\left( 1 +
      \frac{a^{*}_{x-1}}{b^{*}_{x-1}} \right)} \\ & \le
    \frac{1}{\left( 1 + \frac{b^{*}_{x-1}}{a^{*}_x} \right)} +
    \frac{1}{\left( 1 + \frac{a^{*}_{x-1}}{b^{*}_{x-1}} \right)} \\ &
    \le \frac{a^{*}_x}{b^{*}_{x-1}} + \frac{b^{*}_{x-1}}{a^{*}_{x-1}}
    \\ & \le \frac{1}{4(x-1)^2} + \frac{1}{4(x-1)^2} =
    \frac{1}{2(x-1)^2} \;.
\end{align*}
It follows that $\sum_{x = 2}^\infty r_x \le \frac{1}{2}
\sum_{x=2}^\infty \frac{1}{(x-1)^2} < \infty$, so this particular
version of $\Gamma$ is, in fact, trace-class.  However, Proposition
2.6 of \citet{trun:2016} shows that no version of $\Gamma$ is
uniformly ergodic.  Hence, we have a trace-class chain that is not
uniformly ergodic.

\vspace*{8mm}

\noindent{\LARGE \bf Appendix}
\appendix

\section{Proof of Proposition~\ref{prop:tc}}
\label{app:A}

\begin{proof}
By Theorem 2 of QH\&K, it suffices to establish that
\[
\int_{\mathbb{R}^p} k(\beta, \beta) \, d\beta < \infty \;.
\]
Fubini's Theorem implies that
\[
\int_{\mathbb{R}^p} k(\beta, \beta) \, d\beta = \int_{\mathbb{R}_+^n}
\int_{\mathbb{R}^p} \pi(\beta \,|\, w, y) \, \pi(w \,|\, \beta, y) \,
d\beta \, dw \;.
\]
Throughout the proof, let $c_i$, $i = 1, 2, \ldots$, denote finite,
positive constants that do not depend on $\beta$ or $w$.  First,
\begin{align}
  \label{eq:eq1}
    \pi(\beta & \,|\, w, y) \, \pi(w \,|\, \beta, y) \nonumber \\ & =
    c_1 \big \vert \Sigma(w) \big \vert^{-\frac{1}{2}} \mbox{exp}
    \bigg\{ -\frac{1}{2} \Big( \beta^T \big( X^T \Omega(w) X + B^{-1}
    \big) \beta - 2 \beta^T \Big[ X^T \Big( y - \frac{1}{2} 1_n \Big)
      + B^{-1} b \Big] \Big) \bigg\} \nonumber \\ & \;\;\;\; \times
    \mbox{exp} \Big\{ -\frac{1}{2} \mu(w)^T \Sigma(w)^{-1} \mu(w)
    \Big\} \prod_{i=1}^{n} \bigg[ \mbox{cosh} \bigg( \frac{\vert x_i^T
        \beta \vert}{2} \bigg) \mbox{exp} \bigg\{ -\frac{(x_i^T
        \beta)^2}{2} w_i \bigg\} g(w_i) \bigg] \nonumber \\ & \le c_1
    \big \vert \Sigma(w) \big \vert^{-\frac{1}{2}} \mbox{exp} \bigg\{
    -\frac{1}{2} \Big( \beta^T \big( X^T \Omega(w) X + B^{-1} \big)
    \beta - 2 \beta^T \Big[ X^T \Big( y - \frac{1}{2} 1_n \Big) +
      B^{-1} b \Big] \Big) \bigg\} \nonumber \\ & \;\;\;\; \times
    \prod_{i=1}^{n} \bigg[ \mbox{cosh} \bigg( \frac{\vert x_i^T \beta
        \vert}{2} \bigg) \mbox{exp} \bigg\{ -\frac{(x_i^T \beta)^2}{2}
      w_i \bigg\} g(w_i) \bigg] \;.
\end{align}
Now, 
\[
\beta^T X^T 1_n = \sum_{i=1}^n x_i^T \beta
\]
and
\[
\prod_{i=1}^{n} \mbox{cosh} \bigg( \frac{\vert x_i^T \beta \vert}{2}
\bigg) = 2^{-n} \prod_{i=1}^{n} \bigg[ \mbox{exp} \bigg\{ -\frac{x_i^T
    \beta}{2} \bigg\} + \mbox{exp} \bigg\{ \frac{x_i^T \beta}{2}
  \bigg\} \bigg] \;.
\]
Hence,
\begin{equation}
  \label{eq:eq2}
\mbox{exp} \bigg\{ -\frac{1}{2} \beta^T X^T 1_n \bigg\}
\prod_{i=1}^{n} \mbox{cosh} \bigg( \frac{\vert x_i^T \beta \vert}{2}
\bigg) = 2^{-n} \prod_{i=1}^{n} \Big[ 1 + \mbox{exp} \big\{ -x_i^T
  \beta \big\} \Big] \;.
\end{equation}
Also,
\begin{equation}
  \label{eq:eq3}
  \prod_{i=1}^{n} \mbox{exp} \bigg\{ -\frac{(x_i^T \beta)^2}{2} w_i
  \bigg\} = \mbox{exp} \Big\{ -\frac{1}{2} \beta^T X^T \Omega(w) X
  \beta \Big\} \;.
\end{equation}
Combining \eqref{eq:eq1}, \eqref{eq:eq2}, and \eqref{eq:eq3}, we have
\begin{align}
  \label{eq:eq4}
    \pi(\beta \,|\, w, y) \, \pi(w \,|\, \beta, y) \le c_2 \big \vert &
    \Sigma(w) \big \vert^{-\frac{1}{2}} \mbox{exp} \bigg\{
    -\frac{1}{2} \Big( \beta^T \big( 2 X^T \Omega(w) X + B^{-1} \big)
    \beta - 2 \beta^T \Big[ X^T y + B^{-1} b \Big] \Big) \bigg\}
    \nonumber \\ & \times \prod_{i=1}^{n} \Big[ 1 + \mbox{exp} \big\{
      -x_i^T \beta \big\} \Big] g(w_i) \;.
\end{align}
Now, adopting an argument used in \citet{choi:roma:2017}, let $A
\subseteq \mathbb{N}_{n} := \{1, 2, \ldots, n\}$, and define $X_A$ to
be the $n \times p$ matrix whose $i$th row is equal to
\[ 
x_{A, i}^T = \begin{cases} 
              x_i^T & i \in A \\
              0 & \mbox{otherwise} \;.
           \end{cases}
\]
If we let $I_A$ denote a diagonal matrix whose $i$th diagonal element
is 1 if $i \in A$, and 0 otherwise, then we can write $X_A = I_A X$.
Clearly, if $A = \emptyset$, then $X_A$ is a null matrix.  We can see that
\[
\mbox{exp} \big\{ - 1_n^T X_A \beta \big\} = \begin{cases} \mbox{exp}
  \big\{ - \sum_{i \in A} x_i^T \beta \big \} & \mbox{$A$ nonempty}
  \\ 1 & \mbox{otherwise} \;.
   \end{cases}
\]
It then follows that
\begin{align}
  \label{eq:eq5}
  \prod_{i=1}^n \Big[ 1 + \mbox{exp} \big\{ -x_i^T \beta \big\} \Big]
  & = 1 + \sum_{A \subseteq \mathbb{N}_{n} : A \ne \emptyset}
  \mbox{exp} \bigg\{ - \sum_{i \in A} x_i^T \beta \bigg \} \nonumber
  \\ & = \sum_{A \subseteq \mathbb{N}_{n}} \mbox{exp} \big\{ - \beta^T
  X_A^T 1_n \big\} \nonumber \\ & = \sum_{A \subseteq \mathbb{N}_{n}}
  \mbox{exp} \big\{ - \beta^T X^T (I_A 1_n) \big\} \;.
\end{align}
Combining \eqref{eq:eq4} and \eqref{eq:eq5}, and recalling that
$\Sigma(w) = \left( X^T \Omega(w) X + B^{-1} \right)^{-1}$, we have
\begin{align}
  \label{eq:eq6}
    \pi(\beta \,|\, w, y) \, \pi(w \,|\, \beta, y) & \le c_2 \big \vert
    \Sigma(w) \big \vert^{-\frac{1}{2}} \mbox{exp} \Big\{ -\frac{1}{2}
    \beta^T \big( 2 X^T \Omega(w) X + B^{-1} \big) \beta \Big\}
    \nonumber \\ & \;\;\;\; \times \bigg[ \prod_{i=1}^{n} g(w_i)
      \bigg] \sum_{A \subseteq \mathbb{N}_{n}} \mbox{exp} \big\{ -
    \beta^T X^T (I_A 1_n) \big\} \; \mbox{exp} \big\{ \beta^T \big[
      X^T y + B^{-1} b \big] \big\} \nonumber \\ & = c_2 \Bigg[
      \frac{\big \vert X^T \Omega(w) X + B^{-1} \big \vert}{\big \vert
        2 X^T \Omega(w) X + B^{-1} \big \vert} \Bigg]^{\frac{1}{2}}
    \frac{\mbox{exp} \big \{ -\frac{1}{2} \beta^T \big( 2 X^T
      \Omega(w) X + B^{-1} \big) \beta \big \}}{\big \vert \big( 2 X^T
      \Omega(w) X + B^{-1} \big)^{-1} \big \vert^{\frac{1}{2}}}
    \nonumber \\ & \;\;\;\; \times \bigg[ \prod_{i=1}^{n} g(w_i)
      \bigg] \sum_{A \subseteq \mathbb{N}_{n}} \mbox{exp} \big\{
    \beta^T \big[ X^T \big(y - I_A 1_n) + B^{-1} b \big] \big\}
    \nonumber \\ & \le c_2 \frac{\mbox{exp} \big \{ -\frac{1}{2}
      \beta^T \big( 2 X^T \Omega(w) X + B^{-1} \big) \beta \big
      \}}{\big \vert \big( 2 X^T \Omega(w) X + B^{-1} \big)^{-1} \big
      \vert^{\frac{1}{2}}} \nonumber \\ & \;\;\;\; \times \bigg[
      \prod_{i=1}^{n} g(w_i) \bigg] \sum_{A \subseteq \mathbb{N}_{n}}
    \mbox{exp} \big\{ \beta^T \big[ X^T \big(y - I_A 1_n) + B^{-1} b
      \big] \big\} \;,
\end{align}
where the second inequality is a result of the fact that 
\[
\big \vert 2 X^T \Omega(w) X + B^{-1} \big \vert \ge \big \vert X^T
\Omega(w) X + B^{-1} \big \vert + \big \vert X^T \Omega(w) X \big
\vert \ge \vert X^T \Omega(w) X + B^{-1} \big \vert \;,
\]
which follows from the Minkowski determinant inequality \citep[see,
  e.g.,][Theorem 7.8.8]{horn:john:1985}.  Letting $l_A = X^T \big(y -
I_A 1_n) + B^{-1} b$, and using the formula for the moment generating
function of the multivariate normal distribution, we have
\begin{equation}
    \label{eq:eq7}
    \int_{\mathbb{R}^p} e^{\beta^T l_A} \frac{\mbox{exp} \big\{
      -\frac{1}{2} \beta^T \big( 2 X^T \Omega(w) X + B^{-1} \big)
      \beta \big\}}{(2 \pi)^{\frac{n}{2}} \big \vert \big( 2 X^T
      \Omega(w) X + B^{-1} \big)^{-1} \big \vert^{\frac{1}{2}}} \,
    d\beta = \mbox{exp} \big\{l_A^T \big( 2 X^T \Omega(w) X + B^{-1}
    \big)^{-1} l_A/2 \big\} \;.
\end{equation}
Combining \eqref{eq:eq6} and \eqref{eq:eq7} yields
\begin{equation}
  \label{eq:eq8}
    \int_{\mathbb{R}^p} \pi(\beta \,|\, w, y) \, \pi(w \,|\, \beta, y)
    \, d\beta \le c_3 \bigg[ \prod_{i=1}^{n} g(w_i) \bigg] \sum_{A
      \subseteq \mathbb{N}_{n}} \mbox{exp} \big\{l_A^T \big( 2 X^T
    \Omega(w) X + B^{-1} \big)^{-1} l_A/2 \big\} \;.
\end{equation}
Now, since $2 X^T \Omega(w) X + B^{-1} \succeq B^{-1}$, it follows that
$\big( 2 X^T \Omega(w) X + B^{-1} \big)^{-1} \preceq \big( B^{-1}
\big)^{-1} = B$.  Thus,
\[
l_A^T \big( 2 X^T \Omega(w) X + B^{-1} \big)^{-1} l_A \le l_A^T B l_A \;.
\]
Combining this with \eqref{eq:eq8}, we have
\begin{equation}
  \label{eq:eq9}
    \int_{\mathbb{R}^p} \pi(\beta \,|\, w, y) \, \pi(w \,|\, \beta, y)
    \, d\beta \le c_4 \prod_{i=1}^{n} g(w_i) \;,
\end{equation}
and it follows that
\begin{equation*}
\int_{\mathbb{R}_+^n} \int_{\mathbb{R}^p} \pi(\beta \,|\, w, y) \,
\pi(w \,|\, \beta, y) \, d\beta \, dw \le c_4 \int_{\mathbb{R}_+^n}
\prod_{i=1}^{n} g(w_i) \, dw = c_4 \prod_{i=1}^{n} \int_{\mathbb{R}_+}
g(w_i) \, dw_i = c_4 < \infty \;.
\end{equation*}
\end{proof}

\section{Proof of Proposition~\ref{prop:fv}}
\label{app:B}

\begin{proof}
Throughout the proof, let $c_i$, $i = 1, 2, \ldots$, denote finite,
positive constants that do not depend on $\beta$ or $w$.  We begin by
showing that there exists a constant $M \in (0,\infty)$ such that
\begin{equation}
  \label{eq:eq10}
  \frac{\pi(\beta \,|\, w, y)}{h_{\nu}(\beta ; d, C)} \le M \big \vert
  X^T \Omega(w) X + B^{-1} \big \vert^{\frac{1}{2}} \;.
\end{equation}
Letting $z = y - \frac{1}{2} 1_n$, and recalling that $\Sigma(w) =
\left( X^T \Omega(w) X + B^{-1} \right)^{-1}$ and $\Sigma^{-1}(w)
\mu(w) = X^T z + B^{-1} b$, we have
\begin{align*}
  -\frac{1}{2} \big( \beta - \mu(w) \big)^T \Sigma^{-1}(w) \big(
  \beta - \mu(w) \big) & \le -\frac{1}{2} \beta^T \Sigma^{-1}(w) \beta
  + \beta^T \Sigma^{-1}(w) \mu(w) \\ & \le - \frac{1}{2} \beta^T
  B^{-1} \beta + \beta^T B^{-1} \big( B X^Tz + b \big) \\ & =
  - \frac{1}{2} \big( \beta - \big( B X^Tz + b \big) \big)^T B^{-1}
  \big( \beta - \big( B X^Tz + b \big) \big) + c_1 \;.
\end{align*}
It follows that
\begin{align*}
\frac{\pi(\beta \,|\, w, y)}{h_{\nu}(\beta ; d, C)} & = c_2 \big \vert
X^T \Omega(w) X + B^{-1} \big \vert^{\frac{1}{2}} \exp \bigg \{ -
\frac{1}{2} \big( \beta - \mu(w) \big)^T \Sigma^{-1}(w) \big( \beta -
\mu(w) \big) \bigg \} \\ & \hspace*{20mm} \times \bigg( 1 +
\frac{1}{\nu} \big( \beta - d \big)^T C^{-1} \big( \beta - d \big)
\bigg)^{\frac{\nu+p}{2}} \\ & \le c_3 \big \vert X^T \Omega(w) X +
B^{-1} \big \vert^{\frac{1}{2}} \exp \bigg \{ - \frac{1}{2} \big(
\beta - \big( B X^Tz + b \big) \big)^T B^{-1} \big( \beta - \big( B
X^Tz + b \big) \big) \bigg \} \\ & \hspace*{20mm} \times \bigg( 1 +
\frac{1}{\nu} \big( \beta - d \big)^T C^{-1} \big( \beta - d \big)
\bigg)^{\frac{\nu+p}{2}} \\ & \le c_4 \big \vert X^T \Omega(w) X +
B^{-1} \big \vert^{\frac{1}{2}} \;,
\end{align*}
where the last inequality is a consequence of the fact that the
product of the exponential term and the polynomial term converges to 0
as $\beta$ diverges.  Hence, \eqref{eq:eq10} holds.  Now using
\eqref{eq:eq10} and then \eqref{eq:eq9}, we have
\begin{align*}
  \int_{\mathbb{R}_+^n} \int_{\mathbb{R}^p} \frac{ \pi(w \,|\, \beta,
    y) \, \pi^3(\beta \,|\, w, y)}{h_{\nu}^2(\beta ; d, C)} \; d\beta
  \; dw & \le \int_{\mathbb{R}_+^n} M^2 \big \vert X^T \Omega(w) X +
  B^{-1} \big \vert \int_{\mathbb{R}^p} \pi(w \,|\, \beta, y) \,
  \pi(\beta \,|\, w, y) \; d\beta \; dw \\ & \le c_5
  \int_{\mathbb{R}_+^n} \big \vert X^T \Omega(w) X + B^{-1} \big \vert
  \prod_{i=1}^{n} g(w_i) \, dw \;.
\end{align*}
Once expanded, the determinant will result in a finite sum of products
of polynomials of the $w_i$s, and since $\int_{\mathbb{R}_+} u^a g(u)
\, du < \infty$ for all positive integers $a$ \citep[see,
  e.g.,][]{bian:pitm:yor:2001}, the proof is complete.
\end{proof}

\bibliographystyle{ims}
\bibliography{PG}

\end{document}